\documentclass[12pt]{article}

\input{amssym.def}
\input{amssym.tex}

\newcommand{\iy}{\infty}
\newcommand{\bC}{{\bf C}}
\newcommand{\bT}{{\bf T}}
\newcommand{\bZ}{{\bf Z}}

\newcommand{\vsk}{\vspace{1mm}}
\newcommand{\ov}{\overline}
\newcommand{\tic}{\tilde{c}}

\begin{document}

\vspace*{2.5cm}
\begin{center}
{\large\bf ONE MORE PROOF OF THE BORODIN-OKOUNKOV\\[0.5ex]
FORMULA FOR TOEPLITZ DETERMINANTS}
\end{center}

\vspace{4mm}
\begin{center}
{\bf A. B\"ottcher}
\end{center}

\vspace{4mm}
\begin{quote}

\renewcommand{\baselinestretch}{1.0}
\footnotesize
{Recently, Borodin and Okounkov \cite{BO} established a remarkable
identity for Toeplitz determinants. Two other proofs of this identity
were subsequently found by Basor and Widom \cite{BaWi}, who also extended
the formula to the block case. We here give one more proof, also for the
block case. This proof is based on a formula for the inverse of a finite
block Toeplitz matrix obtained in the late seventies by Silbermann and the
author.}
\end{quote}

\vspace{4mm}

\noindent
Given an $N \times N$ matrix function $f$ in $L^\iy$ on the
complex unit circle $\bT$
with Fourier coefficients $\{f_n\}_{n \in \bZ}$, the
Toeplitz operator $T(f)$ and the two Hankel operators $H(f)$
and $H(\tilde{f})$ are defined by the infinite (block) matrices
\[T(f)=(f_{j-k})_{j,k=0}^\iy,\quad H(f)=(f_{j+k+1})_{j,k=0}^\iy,\quad
H(\tilde{f})=(f_{-j-k-1})_{j,k=0}^\iy.\]
These matrices induce bounded operators on $\ell^2(\bZ_+, \bC^N)$.
Let $T_n(f)=(f_{j-k})_{j,k=0}^{n-1}$ and $D_n(f) = \det T_n(f)$.
We denote by $P_n$ the orthogonal projection of $\ell^2(\bZ_+, \bC^N)$
onto the subspace
$\ell^2(\{0,\ldots,n-1\}, \bC^N)$ and we put $Q_n=I-P_n$.
Let finally $K_{2,2}^{1/2,1/2}$ be the Krein algebra of all matrix
functions $f$ in $L^\iy$ on $\bT$ for which $\sum_{n \in \bZ}
(|n|+1)\|f_n\|^2 < \iy$, where $\|\cdot\|$ is any matrix norm.

\vsk
Suppose $a$ is a matrix function in  $K_{2,2}^{1/2,1/2}$ and $a$
admits (left and right Wiener-Hopf) factorizations $a=u_-u_+$ and
$a=v_+v_-$ where $u_+, v_+, \ov{u}_-, \ov{v}_-$ and the inverses
of these matrix functions belong to $K_{2,2}^{1/2,1/2} \cap H^\iy$.
Put $b=v_-u_+^{-1}$ and $c=u_-^{-1}v_+$.  Then the operator
$H(b)H(\tic)$ is of trace class and the operator $I-H(b)H(\tic)$
is invertible. The Borodin-Okounkov formula (\`{a} la Widom) says that
\begin{equation}
D_n(a)=G(a)^n\,\frac{\det(I-Q_nH(b)H(\tic)Q_n)}{\det(I-H(b)H(\tic))}
\label{1}
\end{equation}
for all $n \ge 1$, where $G(a)=\exp \int_0^{2\pi}\log \det a(e^{i
\theta})d\theta$. It is well known that we also have
\begin{equation}
1/\det(I-H(b)H(\tic)) = 1/\det T(b)T(c) = \det T(a)T(a^{-1})
\label{2}
\end{equation}
and that (\ref{2}) equals $\exp \sum_{k=1}^\iy k (\log
a)_k(\log a)_{-k}$ in the scalar case ($N=1)$, where $(\log
a)_j$ is the $j$th Fourier coefficient of $\log a$. In what follows
we abbreviate $H(b)H(\tic)$ to $K$.

\vsk
Our proof is based on the observation that $T_n(a)$ is invertible
if and only if $I-Q_nKQ_n$ is invertible and that in this case
\begin{equation}
T_n^{-1}(a)=T_n(u_+^{-1})\left(I-P_nT(c)Q_n(I-Q_nKQ_n)^{-1}
Q_nT(b)P_n\right)T_n(u_-^{-1}). \label{3}
\end{equation}
If $n$ is large enough then $\|Q_nKQ_n\|<1$, because $K$ is compact
and $Q_n \to 0$ strongly. Hence, for sufficiently large $n$ we can write (\ref{3})
as
\[T_n^{-1}(a)=T_n(u_+^{-1})\left(I-P_nT(c)Q_n\sum_{j=0}^\iy(Q_nKQ_n)^j
Q_nT(b)P_n\right)T_n(u_-^{-1}),\]
and in exactly this form the identity was proved in
\cite[p. 188]{BoSi1} (also see \cite[p. 443]{BoSi2}).

\vsk
Formula (\ref{1}) is immediate from (\ref{3}):
passing  to determinants in (\ref{3}) and taking
into account that
$\det(I+AB)=\det(I+BA)$ and $\det[(I+A)^{-1}]=1/\det(I+A)$, we
get
\begin{eqnarray}
\frac{1}{D_n(a)} & = &
\frac{1}{G(a)^n}\,\det\left(I-Q_nT(b)P_nT(c)Q_n(I-Q_nKQ_n)^{-1}\right)
\nonumber\\
& = & \frac{1}{G(a)^n}\, \frac{\det(I-Q_nKQ_n-Q_nT(b)P_nT(c)Q_n)}
{\det(I-Q_nKQ_n)}, \label{4}
\end{eqnarray}
and since  $I-Q_nKQ_n-Q_nT(b)P_nT(c)Q_n$ equals
\[I-Q_n(T(bc)-T(b)T(c))Q_n-Q_nT(b)P_nT(c)Q_n
= P_n+Q_nT(b)Q_nT(c)Q_n\]
(note that $H(b)H(\tic)=T(bc)-T(b)T(c)$ and $bc =I$), we see that
the numerator in (\ref{4}) is
$\det(P_n+Q_nT(b)Q_nT(c)Q_n)=\det T(b)T(c)$.

\vsk
Here is, for the reader's convenience, a proof of (\ref{3}).
Let $A$ be an invertible operator and let $P$ and $Q$ be complementary
projections. It is well known that the compression $PAP|{\rm Im}\,P$
of $A$ to the range of $P$ is invertible if and only if
$QA^{-1}Q|{\rm Im}\, Q$ is invertible, in which case
\begin{equation}
(PAP)^{-1}P=PA^{-1}P-PA^{-1}Q(QA^{-1}Q)^{-1}QA^{-1}P. \label{5}
\end{equation}
Thus, $T_n(a)$ is invertible if and only if $Q_nT^{-1}(a)Q_n|{\rm Im}\,Q_n$
is invertible, which in turn is equivalent to the invertibility of
\begin{eqnarray*}
& & Q_nT(v_-)Q_n\,Q_nT^{-1}(a)Q_n \,Q_nT(v_+)Q_n |{\rm Im}\,Q_n\\
& & = Q_nT(v_-)Q_nT(u_+^{-1})T(u_-^{-1})Q_nT(v_+)Q_n |{\rm Im}\,Q_n\\
& & =Q_n T(v_-)T(u_+^{-1})T(u_-^{-1})T(v_+)Q_n |{\rm Im}\,Q_n\\
& &=Q_nT(b)T(c)Q_n |{\rm Im}\,Q_n = Q_n|{\rm Im}\,Q_n-Q_nKQ_n|{\rm Im}\,Q_n.
\end{eqnarray*}
Clearly, the last operator is invertible if and only if so is $I-Q_nKQ_n$.
Now suppose that $T_n(a)$ is invertible. From (\ref{5}) and the preceding
computations we obtain
\begin{eqnarray*}
& & T_n^{-1}(a)-P_nT^{-1}(a)P_n\\
& & = P_nT^{-1}(a)Q_n(Q_nT^{-1}(a)Q_n)^{-1}Q_nT^{-1}(a)P_n\\
& & =P_nT(u_+^{-1})P_nT(u_-^{-1})Q_nT(v_+)Q_n(I-Q_nKQ_n)^{-1}
Q_nT(v_-)Q_nT(u_+^{-1})P_nT(u_-)P_n\\
& & = T_n(u_+^{-1})P_nT(c)Q_n(I-Q_nKQ_n)^{-1}Q_nT(b)P_nT_n(u_-^{-1}),
\end{eqnarray*}
and this is (\ref{3}).

\vspace{4mm}
\noindent
\begin{minipage}[t]{7.5cm}
Fakult\" at f\" ur Mathematik \\
Technische Universit\"at Chemnitz \\09107 Chemnitz, Germany \\[0.5ex]
aboettch@mathematik.tu-chemnitz.de
\end{minipage}

\vspace{8mm}
\noindent
MSC 2000: 47B35

\end{document}